\documentclass[a4paper,12 pt]{amsart}
\usepackage{a4,amsbsy,amscd,amssymb,amstext,amsmath,latexsym,oldgerm,enumerate,diagrams}

\makeatletter                                       
\def\l@subsection{\@tocline{2}{0pt}{2.5pc}{2.5pc}{}}
\makeatother                                        

\makeatletter                                                  
\def\chapter{\clearpage\thispagestyle{plain}\global\@topnum\z@ 
\@afterindenttrue \secdef\@chapter\@schapter}
\makeatother

\newtheorem{thmgl} {Theorem}              
\newtheorem{propgl}{Proposition}
\newtheorem{lemgl} {Lemma}

\newtheorem{propnn}{Proposition}
\renewcommand{\thepropnn}{\!\!}

\theoremstyle{definition}

\newtheorem{remnn}{Remark}

\newtheorem{remsnn}{Remarks}

\newcommand{\mf}{\mathfrak}
\newcommand{\mc}{\mathcal}
\newcommand{\mb}{\mathbb}

\newcommand{\nts}{\negthinspace}     
\newcommand{\Nts}{\nts\nts}

\newcommand{\ncd}{\nts\cdot\nts}
\newcommand{\ov}{\overline}

\newcommand{\vink}{^{^{_{_\vee}}}}

\newcommand{\sm}{\setminus}         

\newcommand{\ot}{\otimes}           

\newcommand{\Hom}{{\rm Hom}}        

\newcommand{\Mat}{{\rm Mat}}

\newcommand{\Maxspec}{{\rm Maxspec}}

\newcommand{\Ker}{{\rm Ker}}

\newcommand{\sgn}{{\rm sgn}}
\newcommand{\tr}{{\rm tr}}
\newcommand{\id}{{\rm id}}

\newcommand{\g}{{\mf{g}}}
\newcommand{\h}{{\mf{h}}}

\newcommand{\gl}{\mf{gl}}
\newcommand{\spl}{\mf{sl}}

\newcommand{\GL}{{\rm GL}}
\newcommand{\SL}{{\rm SL}}

\newcommand{\e}{\varepsilon}
\newcommand{\sym}{{\rm sym}}
\begin{document}

\newcommand{\titlesln}{{\Large $\mf{sl}_n$}}
\title{The centre of quantum \titlesln\ at a root of unity}
\author{Rudolf Tange}
\begin{abstract}

It is proved that the centre $Z$ of the simply connected quantised universal enveloping algebra $U_{\e,P}(\spl_n)$, $\e$ a primitive $l$-th root of unity, $l$ an odd integer $>1$, has a rational field of fractions. Furthermore it is proved that if $l$ is a power of an odd prime, $Z$ is a unique factorisation domain. \end{abstract}
\address{Department of Mathematics, The University of Manchester, Oxford Road,
M13 9PL, UK} \email{rudolf.h.tange@stud.man.ac.uk}
\maketitle

\section*{Introduction}

In \cite{DCKP} DeConcini, Kac and Procesi introduced the simply connected quantised universal enveloping algebra $U=U_{\e,P}(\g)$ over $\mb C$ at a primitive $l$-th root of unity $\e$ associated to a simple finite dimensional complex Lie algebra $\g$. The importance of the study of the centre $Z$ of $U$ and its spectrum $\Maxspec(Z)$ is also pointed out in \cite{DCK}.

In this article we consider the following two conjectures concerning the centre $Z$ of $U$ in the case $\g=\spl_n$:
\begin{enumerate}[1.]
\item $Z$ has a rational field of fractions.
\item $Z$ is a unique factorisation domain (UFD).
\end{enumerate}
The same conjectures can be made for the universal enveloping algebra $U(\g)$ of the Lie algebra $\g$ of a reductive group over an algebraically closed field of positive characteristic. In \cite{PrT} these conjectures were proved for $\g=\gl_n$ and for $\g=\spl_n$ under the condition that $n$ is nonzero in the field.

The second conjecture was made by Braun and Hajarnavis in \cite{BHa} for the universal enveloping algebra $U(\g)$ and suggested for $U=U_{\e,P}(\g)$. There it was also proved that $Z$ is locally a UFD. In Section~\ref{s.quf} below, this conjecture is proved for $\spl_n$ under the condition that $l$ is a power of a prime ($\ne 2$). The auxiliary results and step 1 of the proof of Theorem~\ref{thm.quf}, however, hold without extra assumptions on $l$.

The first conjecture was posed as a question by J.\ Alev for the universal enveloping algebra $U(\g)$. It can be considered as a first step towards a proof of a version of the Gelfand-Kirillov conjecture for $U$. Indeed the Gelfand-Kirillov conjecture for $\gl_n$ and $\spl_n$ in positive characteristic\footnote{The {\it Gelfand-Kirillov conjecture} for a Lie algebra $\g$ over $K$ states that the fraction field of $U(\g)$ is isomorphic to a Weyl skew field $D_n(L)$ over a purely transcendental extension $L$ of $K$.}  was proved recently by J.-M.\ Bois in his PhD thesis \cite{Bois} using results in \cite{PrT} on the centres of their universal enveloping algebras (for $\spl_n$ it was required that $n\ne0$ in the field). It should be noted that the Gelfand-Kirillov conjecture for $U(\g)$ in characteristic $0$ (and in positive characteristic) is still open for $\g$ not of type $A$.

As in \cite{PrT}, a certain semi-invariant $d$ for a maximal parabolic subgroup of $\GL_n$ will play an important r\^ole. Later we learned that (a version of) this semi-invariant already appeared before in the literature, see \cite{Dix}. For quantum versions, see \cite{FM1} and \cite{FM2}.

\section{Preliminaries}

In this section we recall some basic results, mostly from \cite{DCKP}, that are needed to prove the main results (Theorems~\ref{thm.qrat} and \ref{thm.quf}) of this article. Short proofs are added in case the results are not explicitly stated in \cite{DCKP}.

\subsection{Elementary definitions}
\ \\
Let $\g$ be a simple finite dimensional Lie algebra over $\mb{C}$ with Cartan subalgebra $\h$, let $\Phi$ be its root system relative to $\h$, let $(\alpha_1,\ldots,\alpha_r)$ be a basis of $\Phi$ and let $(.|.)$ be the symmetric bilinear form on $\h^*$ which is invariant for the Weyl group $W$ and satisfies $(\alpha|\alpha)=2$ for all short roots $\alpha$. Put $d_i=(\alpha_i|\alpha_i)/2$. The root lattice and the weight lattice of $\Phi$ are denoted by resp. $Q$ and $P$. Note that $(.|.)$ is integral on $Q\times P$.

Mostly we will be in the situation where $\g=\spl_n$. In this case $r=n-1$ and all the $d_i$ are equal to 1. We then take $\h$ the subalgebra that consists of the diagonal matrices in $\spl_n$ and we take $\alpha_i=A\mapsto A_{i\,i}-A_{i+1\,i+1}:\h\to\mb{C}$.

Let $l$ be an odd integer $>1$ and coprime to all the $d_i$, let $\e$ be a primitive $l$-th root of unity and let $\Lambda$ be a lattice between $Q$ and $P$. Let $U=U_{\e,\Lambda}(\g)$ be the quantised universal enveloping algebra of $\g$ at the root of unity $\e$ defined in \cite{DCK} and denote the centre of $U$ by $Z$. Since $U$ has no zero divisors (see \cite{DCK} 1.6-1.8), $Z$ is an integral domain. Let $U^+,U^-,U^0$ be the subalgebras of $U$ generated by resp. the $E_i$, the $F_i$ and the $K_\lambda$ with $\lambda\in\Lambda$. Then the multiplication $U^-\ot U^0\ot U^+\to U$ is an isomorphism of vector spaces. We identify $U^0$ with the group algebra $\mb{C}\Lambda$ of $\Lambda$. Note that $W$ acts on $U^0$, since it acts on $\Lambda$. Let $T$ be the complex torus $\Hom(\Lambda,\mb{C}^\times)$. Then $T$ can be identified with $\Maxspec(U^0)=\Hom_{\mb{C}\text{-Alg}}(U^0,\mb{C})$ and for the action of $T$ on $U^0=\mb{C}[T]$ by translation we have $t\cdot K_\lambda=t(\lambda)K_\lambda$.

The {\it braid group} $\mc{B}$ acts on $U$ by automorphisms. See \cite{DCKP} 0.4. The subalgebra $Z_0$ of $U$ is defined as the smallest $\mc{B}$-stable subalgebra containing the elements $K_\lambda^l$, $\lambda\in\Lambda$ and $E_i^l,F_i^l$, $i=1,\ldots,r$. We have $Z_0\subseteq Z$. Put $z_\lambda=K_\lambda^l$ and let $Z_0^0$ be the subalgebra of $Z_0$ spanned by the $z_\lambda$. Then the identification of $U^0$ with $\mb{C}\Lambda$ gives an identification of $Z_0^0$ with $\mb{C}l\Lambda$. If we replace $K_\lambda$ by $z_\lambda$ in foregoing remarks, then we obtain an identification of $T$ with $\Maxspec(Z_0^0)$. Put $Z_0^\pm=Z_0\cap U^\pm$. Then the multiplication $Z_0^-\ot Z_0^0\ot Z_0^+\to Z_0$ is an isomorphism (of algebras). See e.g. \cite{DCK} 3.3.

\subsection{The Harish-Chandra centre $Z_1$ and the quantum restriction theorem}
\ \\
Let $Q\vink$ be the dual root lattice, that is, the $\mb{Z}$-span of the dual root system $\Phi\vink$. We have $Q\vink\cong P^*\hookrightarrow\Lambda^*$. Denote the image of $Q\vink$ under the homomorphism $f\mapsto(\lambda\mapsto (-1)^{f(\lambda)}):\Lambda^*\to T$ by $Q_2\vink$. Then the elements $\ne1$ of $Q_2\vink$ are of order 2 and $U^{0 Q_2\vink}=\mb{C}(\Lambda\cap 2P)$. Since $Q_2\vink$ is $W$-stable, we can form the semi-direct product $\tilde{W}=W\ltimes Q_2\vink$ and then $U^{0 \tilde{W}}=(\mb{C}(\Lambda\cap 2P))^W$.

Let $h':U=U^-\ot U^0\ot U^+\to U^0$ be the linear map taking $x\ot u\ot y$ to $\e_U(x)u\,\e_U(y)$, where $\e_U$ is de counit of $U$. Then $h'$ is a projection of $U$ onto $U^0$. Furthermore $h'(Z_0)=Z_0^0=\mb{C}l\Lambda$ and $h'|_{Z_0}:Z_0\to Z_0^0$ has a similar description as $h'$ and is a homomorphism of algebras. Define the shift automorphism $\gamma$ of $U^{0 Q_2\vink}$ by setting $\gamma(K_\lambda)=\e^{(\rho|\lambda)}K_\lambda$ for $\lambda\in\Lambda\cap 2P$. Here $\rho$ is the half sum of the positive roots. Note that $\gamma=\id$ on $Z_0^{0 Q_2\vink}=\mb{C}l(\Lambda\cap 2P)$. In \cite{DCKP} p 174 and \cite{DCK} \S2, there was constructed an injective homomorphism $\ov{h}:U^{0 \tilde{W}}\to Z$, whose image is denoted by $Z_1$, such that $h'(Z_1)\subseteq U^{0 Q_2\vink}$ and the inverse
$$h:Z_1\stackrel{\sim}{\to}U^{0 \tilde{W}}$$
of $\ov{h}$ is equal to $\gamma^{-1}\circ h'$. Note that $h=h'$ on $Z_0\cap Z_1$ and that $h'|_{Z_1}$ is a homomorphism of algebras. Since $\Ker(h')$ is stable under left and right multiplication by elements of $U^0$ and under multiplication by elements of $Z$, we can conclude that the restriction of $h'$ to the subalgebra generated by $Z_0$ and $Z_1$ is a homomorphism of algebras.

From now on we assume that $\Lambda=P$. Let $G$ be the simply connected almost simple complex algebraic group with Lie algebra $\g$ and let $T$ be a maximal torus of $G$. We identify $\Phi$ and $W$ with the root system and the Weyl group of $G$ relative to $T$. Note that the character group $X(T)$ of $T$ is equal to $P$. In case $\g=\spl_n$ we take $T$ the subgroup of diagonal matrices in $\SL_n$.

\subsection{Generators for $\mb{C}[G]^G$ and $Z_1$}\label{ss.groupinvgen}
\ \\
We denote the fundamental weights corresponding to the basis $(\alpha_1,\ldots,\alpha_r)$ by $\varpi_1,\ldots,\varpi_r$. As is well known, they form a basis of $P$. Let $\mb{C}[G]$ be the algebra of regular functions on $G$. Then the restriction homomorphism $\mb{C}[G]\to\mb{C}[T]=\mb{C}P$ induces an isomorphism $\mb{C}[G]^G\stackrel{\sim}{\to}\mb{C}[T]^W=(\mb{C}P)^W$, see \cite{St1} \S6. For $\lambda\in P$ denote the basis element of $\mb{C}P$ corresponding to $\lambda$ by $e(\lambda)$, denote the $W$-orbit of $\lambda$ by $W\ncd\lambda$ and put $\sym(\lambda)=\sum_{\mu\in W\ncd\lambda}e(\mu)$. Then the ${\rm sym}(\varpi_i)$, $i=1,\ldots,r$ are algebraically independent generators of $(\mb{C}P)^W$. See \cite{Bou1} no. VI.3.4, Thm. 1.

For a field $K$, we denote the vector space of all $n\times n$ matrices over $K$ by $\Mat_n=\Mat_n(K)$. Now assume that $K={\mb C}$. In this section we denote the restriction to $\SL_n$ of the standard coordinate functionals on $\Mat_n$ by $\xi_{ij}$, $1\le i,j\le n$. Furthermore, for $i\in\{1,\ldots,n-1\}$, $s_i\in\mb{C}[\SL_n]$ is defined by $s_i(A)=\tr(\wedge^iA)$, where $\wedge^iA$ denotes the $i$-th exterior power of $A$ and $\tr$ denotes the trace. Then $\varpi_i=(\xi_{11}\cdots\xi_{ii})|_T$ and therefore ${\rm sym}(\varpi_i)=s_i|_T$~(*), the $i$-th elementary symmetric function in the $\xi_{jj}|_T$. See \cite{PrT} 2.4.

In the general case we use the restriction theorem for $\mb{C}[G]$ and define $s_i\in\mb{C}[G]^G$ by (*). So then $s_1,\ldots,s_r$ are algebraically independent generators of $\mb{C}[G]^G$.

Identifying $U^0$ and $\mb{C}P$, we have $U^{0 \tilde{W}}=(\mb{C}2P)^W$. Put $u_i=\ov{h}(\sym(2\varpi_i))$. Then $h(u_i)=\sym(2\varpi_i)$ and $u_1,\ldots,u_r$ are algebraically independent generators of $Z_1$.

\subsection{The cover $\pi$ and the intersection $Z_0\cap Z_1$}
\ \\
Let $\Phi^+$ be the set of positive roots corresponding to the basis $(\alpha_1,\ldots,\alpha_r)$ of $\Phi$ and let $U_+$ resp $U_-$ be the maximal unipotent subgroup of $G$ corresponding to $\Phi^+$ resp. $-\Phi^+$. If $\g=\spl_n$, then $U_+$ and $U_-$ consist of the upper resp. lower triangular matrices in $\SL_n$ with ones on the diagonal. Put $\mc{O}=U_-TU_+$. Then $\mc{O}$ is a nonempty open and therefore dense subset of $G$. Furthermore, the group multiplication defines an isomorphism $U_-\times T\times U_+\stackrel{\sim}{\to}\mc{O}$ of varieties. Put $\Omega=\Maxspec(Z_0)$.

In \cite{DCK} (3.4-3.6) there was constructed a group $\tilde{G}$ of automorphisms of $\hat{U}=\hat{Z_0}\ot_{Z_0}U$, where $\hat{Z}_0$ denotes the algebra of holomorphic functions on the complex analytic variety $\Omega$. The group $\tilde{G}$ leaves $\hat{Z_0}$ and $\hat{Z}=\hat{Z_0}\ot_{Z_0}Z$ stable. In particular it acts by automorphisms on the complex analytic variety $\Omega$. In \cite{DCKP} this action is called the "quantum coadjoint action".

In \cite{DCKP} \S4 there was constructed an unramified cover $\pi:\Omega\to \mc{O}$ of degree $2^r$. I give a short description of the construction of $\pi$. Put $\Omega^\pm=\Maxspec(Z_0^\pm)$. Then we have $\Omega=\Omega^-\times T\times\Omega^+$. Now $Z:\Omega\to T$ is defined as the projection on $T$, $X:\Omega\to U_+$ and $Y:\Omega\to U_-$ as the projection on $\Omega^\pm$ followed by some isomorphism $\Omega^\pm\stackrel{\sim}{\to}U_\pm$ and $\pi$ is defined as $YZ^2X$ (multiplication in $G$).\footnote{In \cite{DCKP} $Z^2$ is denoted by $Z$. The notation here comes from \cite{DCP}. The centre of $U$ is denoted by the same letter, but this will cause no confusion.} This means: $\pi(x)=Y(x)Z(x)^2X(x)$.

The following proposition says something about how $\tilde{G}$ and $\pi$ are related to the "Harish-Chandra centre" $Z_1$ and the conjugation action of $G$ on $\mb{C}[G]$. For more precise statements see 5.4, 5.5 and \S6 in \cite{DCKP}.
\begin{thmgl}[\cite{DCKP} Prop 6.3, Thm 6.7]\label{thm.coverHCcentre}
Consider $\pi$ as a morphism to $G$. Then the comorphism $\pi^{co}:\mb{C}[G]\to Z_0$ is injective and the following holds:
\begin{enumerate}[(i)]
\item $Z^{\tilde{G}}=Z_1$.\ \footnotemark
\item $\pi^{co}$ induces an isomorphism $\mb{C}[G]^G\stackrel{\sim}{\to}Z_0^{\tilde{G}}=Z_0\cap Z_1$.
\item The monomorphism $(\mb{C}P)^W\stackrel{\sim}{\to}(\mb{C}P)^W$ obtained by combining the isomorphism in (ii) with the restriction homomorphism $\mb{C}[G]\to\mb{C}[T]={\mb C}P$ and $h:Z_1\to U^0={\mb C}P$, is given by $x\mapsto 2lx:P\to P$. In particular $h(Z_0\cap Z_1)=(\mb{C}2lP)^W$.
\end{enumerate}
\end{thmgl}
\footnotetext{$\tilde{G}$ is a group of automorphisms of the algebra $\hat{U}$ and does not leave $Z$ stable. However, $S^{\tilde{G}}$ can be defined in the obvious way for every subset $S$ of $\hat{U}$.} I will give the proof of (iii). If we identify $Z_0^0$ with ${\mb C}[T]$, then the homomorphism $h'|_{Z_0}:Z_0\to Z_0^0$ is the comorphism of a natural embedding $T\hookrightarrow\Omega$. Now we have a commutative diagram
\begin{diagram}[small,nohug]
G      &\lTo^{\pi}         &\Omega\\
\uIntoB&                   &\uIntoB\\
T      &\lTo^{t\mapsto t^2}&T\\
\end{diagram}
Expressed in terms of the comorphisms this reads: $(x\mapsto 2x)\circ{\rm res}_{G,T}={\rm res}_{\Omega,T}\circ\pi^{co}$, where ${\rm res}_{G,T}$ and ${\rm res}_{\Omega,T}$ are the restrictions to $T$ and the comorphism of the morphism between the tori is denoted by its restrictions to the character groups. Now we identify $U^0$ with ${\mb C}[T]$. Composing both sides on the left with $x\mapsto lx$ and using $(x\mapsto lx)\circ{\rm res}_{\Omega,T}=h'|_{Z_0}:Z_0\to U^0=\mb{C}P$ we obtain $(x\mapsto 2lx)\circ{\rm res}_{G,T}=h'\circ\pi^{co}$. If we restrict both sides of this equality to $\mb{C}[G]^G$, then we can replace $h'$ by $h$ and we obtain the assertion.\\

\subsection{$Z_0$ and $Z_1$ generate $Z$}\label{ss.Z0andZ1}
\begin{thmgl}[\cite{DCKP} Proposition 6.4, Theorem 6.4]\label{thm.Z0andZ1}
Let $u_1,\ldots,u_r$ be the elements of $Z_1$ defined in Subsection~\ref{ss.groupinvgen}. Then the following holds:
\begin{enumerate}[(i)]
\item The multiplication $Z_1\ot_{Z_0\cap Z_1}Z_0\to Z$ is an isomorphism of algebras.
\item $Z$ is a free $Z_0$-module of rank $l^r$ with the restricted monomials $u_1^{k_1}\cdots u_r^{k_r}$, $0\leq k_i<l$ as a basis.
\end{enumerate}
\end{thmgl}
I give a proof of (ii). In \cite{DCKP} Prop.~6.4 it is proved that $(\mb{C}P)^W$ is a free $(\mb{C}lP)^W$-module of rank $l^r$ with the restricted monomials (exponents $<l$) in the ${\rm sym}(\varpi_i)$ as a basis. The same holds then of course for $(\mb{C}2P)^W$, $(\mb{C}2lP)^W$ and the ${\rm sym}(2\varpi_i)$. But then the same holds for $Z_1$, $Z_0\cap Z_1$ and the $u_i$ by (iii) of Theorem~\ref{thm.coverHCcentre}. So the result follows from (i).

Recall that $\Omega=\Omega^-\times T\times\Omega^+$ and that $\Omega^\pm\cong U_\pm$. So $Z_0$ is a polynomial algebra in $\dim(\g)$ variables with $r$ variables inverted. In particular it's Krull dimension (which coincides with the transcendence degree of its field of fractions) is $\dim(\g)$. The same holds then for $Z$, since it is a finitely generated $Z_0$-module.

Let $Z_0'$ be a subalgebra of $Z_0$ containing $Z_1\cap Z_0$. Then the multiplication $Z_1\ot_{Z_0\cap Z_1}Z_0'\to Z_0'Z_1$ is an isomorphism of algebras by the above theorem. This gives us a way to determine generators and relations for $Z_0'Z_1$: Let $s_1,\ldots,s_r$ be the generators of ${\mb C}[G]^G$ defined in Subsection~\ref{ss.groupinvgen}. Then $\pi^{co}(s_1),\ldots,\pi^{co}(s_r)$ are generators of $Z_0\cap Z_1=Z_0'\cap Z_1$ by Theorem~\ref{thm.coverHCcentre}(ii). Now assume that we have generators and relations for $Z_0'$. We use for $Z_1$ the generators $u_1,\ldots,u_r$ defined in Subsection~\ref{ss.groupinvgen}. For each $i\in\{1,\ldots,r\}$ we can express $\pi^{co}(s_i)$ as a polynomial $g_i$ in the generators of $Z_0'$ and as a polynomial $f_i$ in the $u_j$. Then the generators and relations for $Z_0'$ together with the $u_i$ and the relations $f_i=g_i$ form a presentation of $Z_0'Z_1$.\footnote{This method was also used by Krylyuk in \cite{Kry} to determine generators and relations for the centre of the universal enveloping algebra $U(\g)$ of $\g$. Our homomorphism $\pi^{co}:\mb{C}[G]\to Z_0$ plays the r\^ole of Krylyuk's $G$-equivariant isomorphism $\eta:S(\g)^{(1)}\to Z_p$, where we use the notation of \cite{PrT}.}

The $f_i$ can be determined as follows. Write $\sym(l\varpi_i)$ as a polynomial $f_i$ in the $\sym(\varpi_j)$. Then $\sym(2l\varpi_i)$ is the same polynomial in the $\sym(2\varpi_j)$ and $\pi^{co}(s_i)=f_i(u_1,\ldots,u_r)$ by Theorem~\ref{thm.coverHCcentre}(iii).

Note that $\pi^{co}(\mb{C}[\mc{O}])=Z_0^-\mb{C}(2lP)Z_0^+$ and that $Z_0=\pi^{co}(\mb{C}[\mc{O}])[z_{\varpi_1},\ldots,z_{\varpi_r}]$.

Now assume that $G=\SL_n$. For $f\in\mb{C}[\SL_n]$ denote $f\circ\pi$ by $\tilde{f}$ and put $\tilde{Z_0}=\pi^{co}(\mb{C}[\SL_n])$. Then $\tilde{Z_0}$ is generated by the $\tilde{\xi}_{ij}$; it is a copy of $\mb{C}[\SL_n]$ in $Z_0$. Now $\mc{O}$ consists of the matrices $A\in\SL_n$ that have an LU-decomposition (without row permutations), that is, whose principal minors $\Delta_1(A),\ldots,\Delta_{n-1}(A)$ are nonzero. So $\mb{C}[\mc{O}]= \mb{C}[\SL_n][\Delta_1^{-1},\ldots,\Delta_{n-1}^{-1}]$, $\pi^{co}(\mb{C}[\mc{O}])= \tilde{Z_0}[\tilde{\Delta}^{-1}_1,\ldots,\tilde{\Delta}^{-1}_{n-1}]$ and
$$Z_0=\tilde{Z_0}[z_{\varpi_1},\ldots,z_{\varpi_{n-1}}][\tilde{\Delta}^{-1}_1,\ldots, \tilde{\Delta}^{-1}_{n-1}].$$
Let ${\rm pr}_{\mc{O},T}$ be the projection of $\mc{O}$ on $T$. An easy computation shows that $\Delta_i|_\mc{O}=(\xi_{11}\cdots\xi_{ii})\circ{\rm pr}_{\mc{O},T}=\varpi_i\circ{\rm pr}_{\mc{O},T}$ for $i=1,\ldots,n-1$.\footnote{For two $n\times n$ matrices $A$ and $B$ we have $\wedge^k(AB)=\wedge^k(A)\wedge^k(B)$. From this it follows that if either $A$ is lower triangular or $B$ is upper triangular, then $\Delta_k(AB)=\Delta_k(A)\Delta_k(B)$.} So $\tilde{\Delta}_i=\varpi_i\circ {\rm pr}_{\mc{O},T}\circ\pi=\varpi_i\circ(t\mapsto t^2)\circ {\rm pr}_{\Omega,T}=2\varpi_i\circ {\rm pr}_{\Omega,T}=z_{\varpi_i}^2$. In Subsection~\ref{ss.presentation} we will determine generators and relations for $Z_0'Z_1$, where $Z_0'=\tilde{Z_0}[z_{\varpi_1},\ldots,z_{\varpi_{n-1}}]$ using the method mentioned above.


\section{Rationality}
We use the notation of Section 1 with the following modifications. The functions $\xi_{ij}$, $1\le i,j\le n$, now denote the standard coordinate functionals on $\Mat_n$ and for $i\in\{1,\ldots,n\}$, $s_i\in K[\Mat_n]$ is defined by $s_i(A)=\tr(\wedge^iA)$ for $A\in\Mat_n$. Then $\det(x\id-A)=$ $x^n+\sum_{i=1}^n (-1)^i s_i(A)x^{n-i}$. This notation is in accordance with \cite{PrT}.

For $f\in\mb{C}[\Mat_n]$ we denote its restriction to $\SL_n$ by $f'$ and we denote $\pi^{co}(f')$ by $\tilde{f}$. So now $s_1',\ldots,s_{n-1}'$ and $\xi_{ij}'$ are the functions $s_1,\ldots,s_{n-1}$ and $\xi_{ij}$ of Subsection~\ref{ss.groupinvgen} and the $\tilde{\xi_{ij}}$ are the same.

To prove theorem below we need to look at the expressions of the functions $s_i$ in terms of the $\xi_{ij}$. We use that those equations are linear in $\xi_{1n},\xi_{2n},\ldots,\xi_{nn}$. The treatment is completely analogous to that in \cite{PrT} 4.1 (we use the same symbols $R$, $M$, $d$ and $x_{\bf a}$) to which we refer for more explanation. Let $R$ be the $\mb{Z}$-subalgebra of $\mb{C}[\Mat_n]$ generated by all $\xi_{ij}$ with $j\ne n$.

Let $\partial_{ij}$ denote differentiation with respect to the variable $\xi_{ij}$ and set
$$M= \left[\begin{array}{cccc}
\partial_{1n} (s_1)&\partial_{2n}(s_1)&\ldots& \partial_{nn}(s_1)\\
\partial_{1n} (s_2)&\partial_{2n}(s_2)&\ldots& \partial_{nn}(s_2)\\
\vdots&\vdots& &\vdots \\
\partial_{1n} (s_n)&\partial_{2n}(s_n)&\ldots& \partial_{nn}(s_n)
\end{array}\right],\quad
{\bf c}=\left[\begin{array}{c}
\xi_{1n}\\
\xi_{2n}\\
\vdots \\
\xi_{nn}
\end{array}\right],\quad
{\bf s}=\left[\begin{array}{c}
s_1\\
s_2\\
\vdots \\
s_n
\end{array}\right].
$$
Then the matrix $M$ has entries in $R$ and the following vector equation holds:
\begin{equation}\label{eq.gln2}
M\cdot {\bf c}\,=\,{\bf s} + {\bf r},  \quad\,\mbox{where} \  \ {\mathbf r}\in R^n.
\end{equation}
We denote the determinant of $M$ by $d$. For ${\bf a}=(a_1,\ldots,a_n)\in K^n$ we set
$$x_{\bf a}\,=\,
\left[\begin{array}{ccccc}
0&\cdots&0&0&a_n\\
1&\cdots&0&0& a_{n-1}\\
\vdots&\ddots&\vdots&\vdots&\vdots\\
0&\cdots&1&0&a_2 \\
0&\cdots&0&1&a_1\end{array}\right].$$ Then the minimal polynomial of $x_{\bf a}$ equals $x^n-\sum_{i=1}^n a_ix^{n-i}$, $\det(x_{\bf a})=(-1)^{n-1}a_n$ and $d(x_{\bf{a}})=1$ (Compare Lemma~3 in \cite{PrT}).

\begin{thmgl}\label{thm.qrat}
$Z$ has a rational field of fractions.
\end{thmgl}

\begin{proof}
Denote the field of fractions of $Z$ by $Q(Z)$. From Subsection~\ref{ss.Z0andZ1} it is clear that $Q(Z)$ has transcendence degree $\dim(\spl_n)=n^2-1$ over $\mb{C}$ and that it is generated as a field by the $n^2+2(n-1)$ variables $\tilde{\xi}_{ij}$, $u_1,\ldots,u_{n-1}$ and $z_{\varpi_1},\ldots,z_{\varpi_{n-1}}$. To prove the assertion we will show that $Q(Z)$ is generated by the $n^2-1$ elements $\tilde{\xi}_{ij}$, $i\ne j$, $j\ne n$, $u_1,\ldots,u_{n-1}$ and $z_{\varpi_1},\ldots,z_{\varpi_{n-1}}$. We will first eliminate the $n$ generators $\tilde{\xi}_{1n},\ldots,\tilde{\xi}_{nn}$ and then the $n-1$ generators $\tilde{\xi}_{11},\ldots,\tilde{\xi}_{n-1,n-1}$.

Applying the homomorphism $f\mapsto\tilde{f}=\pi^{co}\circ (f\mapsto f'):\mb{C}[\Mat_n]\to Z_0$ to both sides of \eqref{eq.gln2} we obtain the following equations in the $\tilde{\xi}_{ij}$ and $\tilde{s}_1,\ldots,\tilde{s}_{n-1}$
\begin{equation}\label{eq.Ze}
\tilde{M}\cdot \tilde{{\bf c}}\,=\,\tilde{{\bf s}} + \tilde{{\bf r}},\text{\quad where\ }\tilde{{\bf r}}\in \tilde{R}^n.
\end{equation}
Here $\tilde{M},\tilde{{\bf c}},\tilde{{\bf s}},\tilde{{\bf r}}$ have the obvious meaning, except that we put the last component of $\tilde{{\bf s}}$ and $\tilde{{\bf r}}$ equal to 0 resp. 1, and $\tilde{R}$ is the $\mb{Z}$-subalgebra of $Z_0$ generated by all $\tilde{\xi}_{ij}$ with $j\ne n$. Choosing $\bf{a}$ such that $a_n=(-1)^{n-1}$ we have $x_{\bf a}\in\SL_n$. Since $d(x_{\bf a})=1$, we have $d'\ne 0$ and therefore $\det(\tilde{M})=\tilde{d}\ne0$. Furthermore, for $i=1,\ldots,n-1$, $(\tilde{{\bf s}})_i=\tilde{s}_i\in Z_0\cap Z_1$ and $Z_1$ is generated by $u_1,\ldots,u_{n-1}$. It follows that $\tilde{\xi}_{1n},\ldots,\tilde{\xi}_{nn}$ are in the subfield of $Q(Z)$ generated by the $\tilde{\xi}_{ij}$ with $j\ne n$ and $u_1,\ldots,u_{n-1}$.

Now we will eliminate the generators $\tilde{\xi}_{11},\ldots,\tilde{\xi}_{n-1,n-1}$. We have $$z_{\varpi_1}^2=\tilde{\Delta}_1=\tilde{\xi}_{11}$$ and for $k=2,\ldots,n-1$ we have, by the Laplace expansion rule,
$$z_{\varpi_k}^2= \tilde{\Delta}_k= \tilde{\xi}_{kk}\tilde{\Delta}_{k-1}+t_k= \tilde{\xi}_{kk}z_{\varpi_{k-1}}^2+t_k,$$
where $t_k$ is in the $\mb{Z}$-subalgebra of $Z$ generated by the $\tilde{\xi}_{ij}$ with $i,j\le k$ and $(i,j)\ne (k,k)$. It follows by induction $k$ that for $k=1,\ldots,n-1$, $\tilde{\xi}_{11},\ldots,\tilde{\xi}_{kk}$ are in the subfield of $Q(Z)$ generated by the $z_{\varpi_i}$ with $i\le k$ and the $\tilde{\xi}_{ij}$ with $i,j\le k$ and $i\ne j$.
\end{proof}

\section{Unique Factorisation}\label{s.quf}
Recall that {\it Nagata's lemma} asserts the following: If $x$ is a prime element of a Noetherian integral domain $S$ such that $S[x^{-1}]$ is a UFD, then $S$ is a UFD. See \cite{Eis} Lemma 19.20. In Theorem~\ref{thm.quf} we will see that, by Nagata's lemma, it suffices to show that the algebra $Z/(\tilde{d})$ is an integral domain in order to prove that $Z$ is a UFD. To prove this we will show by induction that the two sequences of algebras (to be defined later): $$K[\SL_n]/(d')\cong\ov{A}(K)=\ov{B}_{0,0}(K)\subseteq\ov{B}_{0,1}(K)\subseteq\cdots\subseteq \ov{B}_{0,n-1}(K)=\ov{B}_0(K)$$ in characteristic $p$ and
$$\ov{B}_0({\mb C})\subseteq\ov{B}_1({\mb C})\subseteq\cdots\subseteq\ov{B}_{n-1}({\mb C})=\ov{B}({\mb C})$$
consist of integral domains. Lemma's~\ref{lem.sln} and \ref{lem.slnd} are, among other things, needed to show that $\ov{A}(K)\cong K[\SL_n]/(d')$ is an integral domain. Lemma's~\ref{lem.leadmondetd} and \ref{lem.leadmonfi} are needed to obtain bases over $\mb Z$ (see Proposition~\ref{prop.presentation(d)}), which, in turn, is needed to pass to fields of positive characteristic and to apply mod $p$ reduction (see Lemma~\ref{lem.modp}).

\subsection{The case $n=2$}\label{ss.nis2}
\ \\
In this subsection we show that the centre of $U_{\e,P}(\spl_2)$ is always a UFD, without any extra assumptions on $l$. The standard generators for $U=U_{\e,P}(\spl_2)$ are $E,F,K_\varpi$ and $K_\varpi^{-1}$. Put $K=K_\alpha=K_\varpi^2$, $z_1=z_\varpi=K_\varpi^l$, $z=z_\alpha=z_1^2=K^l$. Furthermore, following \cite{DCKP} 3.1, we put $c=(\e-\e^{-1})^l$, $x=-cz^{-1}E^l$, $y=cF^l$. Then $x,y$ and $z_1$ are algebraically independent over $\mb C$ and $Z_0={\mb C}[x,y,z_1][z_1^{-1}]$ (see \cite{DCKP} \S 3).

We have $U^0=\mb{C}[K_\varpi,K_\varpi^{-1}]$ and $U^{0\tilde{W}}=\mb{C}[K,K^{-1}]^W=\mb{C}[K+K^{-1}]$. Identifying $U^0$ and $\mb{C}P$, we have $\sym(2\varpi)=K+K^{-1}$ and $\sym(2l\varpi)=z+z^{-1}$. Put $u=\ov{h}(\sym(2\varpi))$. By the restriction theorem for $U$, $Z_1$ is a polynomial algebra in $u$. Denote the trace map on $\Mat_2$ by $\tr$. Then $\tr|_T=\sym(\varpi)$. By the restriction theorem for ${\mb C}[G]$ and Theorem~\ref{thm.coverHCcentre}(ii), $\tilde{\tr}$ generates $Z_0\cap Z_1$. Furthermore $\tilde{\tr}=\ov{h}(z+z^{-1})$, by Theorem~\ref{thm.coverHCcentre}(iii). Let $f\in\mb{C}[u]$ be the polynomial with $z+z^{-1}=f(K+K^{-1})$. Then $\tilde{\tr}=f(u)$. From the formula's in \cite{DCKP} 5.2 it follows that $\tilde{\tr}=-zxy+z+z^{-1}$.

By the construction from Subsection~\ref{ss.Z0andZ1} (we take $Z_0'=Z_0$), $Z$ is isomorphic to the quotient of the localised polynomial algebra $\mb{C}[x,y,z_1,u][z_1^{-1}]$ by the ideal generated by $-z_1^2xy+z_1^2+z_1^{-2}-f(u)$. Clearly $x,u$ and $z_1$ generate the field of fractions of $Z$. In particular they are algebraically independent. So $Z[x^{-1}]$ is isomorphic to the localised polynomial algebra $\mb{C}[x,z_1,u][z_1^{-1},x^{-1}]$ and therefore a UFD. By Nagata's lemma it suffices to show that $x$ is a prime element in $Z$. But $Z/(x)$ is isomorphic to the quotient of $\mb{C}[y,z_1,u][z_1^{-1}]$ by the ideal generated by $z_1^2+z_1^{-2}-f(u)$. This ideal is also generated by $z_1^4-f(u)z_1^2+1$. So it suffices to show that $z_1^4-f(u)z_1^2+1$ is irreducible in $\mb{C}[y,z_1,u][z_1^{-1}]$. From the fact that $f$ is of odd degree $l>0$ (see e.g. Lemma~\ref{lem.leadmonfi} below), one easily deduces that $z_1^4-f(u)z_1^2+1$ is irreducible in $\mb{C}[z_1,u]$ %
and therefore also in $\mb{C}[y,z_1,u]$. Clearly $z_1^4-f(u)z_1^2+1$ is not invertible in $\mb{C}[y,z_1,u][z_1^{-1}]$, so it is also irreducible in this ring.

\subsection{$\SL_n$ and the function $d$}\label{ss.slnd}
\ \\
Part (i) of the next lemma is needed for the proof of Lemma~\ref{lem.sln} and part (ii) is needed for the proof of Theorem~\ref{thm.quf}. The Jacobian matrices below consist of the partial derivatives of the functions in question with respect to the variables $\xi_{ij}$.
\begin{lemgl}\label{lem.jac}
\begin{enumerate}[(i)]
\item There exists a matrix $A\in\SL_n({\mb Z})$ such that $\Delta_{n-1}(A)=0$ and such that some second order minor of the Jacobian matrix of $(\det,\Delta_{n-1})$ is $\pm1$ at $A$.
\item If $n\ge3$, then there exists a matrix $A\in\SL_n({\mb Z})$ such that $d(A)=0$ and such that some $2n$-th order minor of the Jacobian matrix of $(s_1,\ldots,s_n,d,\Delta_1,\ldots,$\\ $\Delta_{n-1})$ is $\pm1$ at $A$.
\end{enumerate}
\end{lemgl}

\begin{proof}
The computations below are very similar to those in \cite{PrT} Section~6.
We denote by $\mc{X}$ the $n\times n$ matrix $(\xi_{ij})$ and for an $n\times n$ matrix $B=(b_{ij})$ and $\Lambda_1 ,\Lambda_2\subseteq\{1,\ldots,n\}$ we denote by $B_{\Lambda_1,\,\Lambda_2}$ the matrix $(b_{ij})_{i\in\Lambda_1,j\in\Lambda_2}$, where the indices are taken in the natural order.

In the computations below we will use the following two facts:\\
For $\Lambda_1,\Lambda_2\subseteq\{1,\ldots,n\}$ with $|\Lambda_1|=|\Lambda_2|$ we have
$$\partial_{ij}\big(\det(\mc{X}_{\Lambda_1,\,\Lambda_2})\big)=
\begin{cases}
(-1)^{n_1(i)+n_2(j)}\det(\mc{X}_{\Lambda_1\sm\{i\},\,\,\Lambda_2\sm\{j\}})&
\mbox{when}\ \, (i,j)\in (\Lambda_1\times\Lambda_2),\\
0&\mbox{when}\ \,(i,j)\not\in (\Lambda_1\times \Lambda_2),
\end{cases}$$
where $n_1(i)$ denotes the position in which $i$ occurs in $\Lambda_1$ and similarly for $n_2(j)$.\\
For $k\le n$ we have $s_k=\sum_{\Lambda}\det(\mc{X}_{\Lambda,\,\Lambda})$ where the sum ranges over all $k$-subsets $\Lambda$ of $\{1,\ldots,n\}$.\\
(i).\ We let $A$ be the following $n\times n$-matrix:
$$A=\left[\begin{array}{ccccc}
0&0&\cdots&0&-1\\
0&1&\cdots&0&0\\
\vdots&\vdots&\ddots&\vdots&\vdots\\
0&0&\cdots&1&0\\
1&0&\cdots&0&0\end{array}\right].$$ Clearly $\det(A)=1$ and $\Delta_{n-1}(A)=0$. From the above two facts it is easy to deduce that $\left[\begin{array}{cc}
\partial_{1\,n}\det&\partial_{1\,1}\det\\
\partial_{1\,n}\Delta_{n-1}&\partial_{1\,1}\Delta_{n-1}
\end{array}\right]$ is equal to
$\left[\begin{array}{cc}
\pm1&0\\
0&\pm1
\end{array}\right]$ at $A$.\\
(ii).\ Put $\alpha=\big((1\,1),(2\,2),(2\,3),\ldots,(2\,n-1),(n\,n),(n-1\,n),\ldots,(2\,n),(2\,1),(1\,2)\big)$, and let $\alpha_i$ denote the $i$-th component of $\alpha$. We let $A$ be the following $n\times n$-matrix:
$$A=\left[\begin{array}{cccccc}
1&0&0&\cdots&0&\Nts(-1)^n\Nts\\
0&1&0&\cdots&0&0\\
1&1&0&\cdots&0&0\\
0&0&1&\cdots&0&0\\
\vdots&\vdots&\vdots&\ddots&\vdots&\vdots\\
0&0&0&\cdots&1&0\end{array}\right].$$
The columns of the Jacobian matrix are indexed by the pairs $(i,j)$ with $1\le i,j\le n$. Let $M_\alpha$ be the $2n$-square submatrix of the Jacobian matrix consisting of the columns with indices from $\alpha$. By permuting in $A$ the first row to the last position and interchanging the first two columns, we see that $\det(A)=1$. We will show that $d(A)=0$ and that the minor $d_{\alpha}:=\det(M_\alpha)$ of the Jacobian matrix is nonzero at $A$.

First we consider the $\Delta_k$, $k\in\{1,\ldots,n-1\}$. By inspecting the matrix $A$ and using the fact that $\partial_{ij}\Delta_k=0$ if $i>k$ or $j>k$, we deduce the following facts: $(\partial_{2\,i}\Delta_k)(A)=
\begin{cases}
\pm1&\text{if $i=k$,}\\
0&\text{if $i>k$,}
\end{cases}$\quad for $i,k\in\{1,\ldots,n-1\}\ i\ne1$,\quad $(\partial_{1\,1}\Delta_1)(A)=1$,\\
$(\partial_{1\,2}\Delta_k)(A)=(\partial_{2\,1}\Delta_k)(A)=0$ for all $k\in\{1,\ldots,n-1\}$ and\\
$(\partial_{i\,n}\Delta_k)(A)=0$ for all $k\in\{1,\ldots,n-1\}$ and all $i\in\{1,\ldots,n\}$.

Now we consider the $s_k$. Let $i\in\{i,\ldots,n\}$ and let $\Lambda\subseteq\{1,\ldots,n\}$. Assume that $\partial_{i n}\big(\det(\mc{X}_{\Lambda,\,\Lambda})\big)$ is nonzero at $A$. Then we have:
\begin{enumerate}[$\bullet$]
\item $i,n\in\Lambda$;
\item $j\in\Lambda\Rightarrow j-1\in\Lambda$ for all $j$ with $4\leq j\leq n$ and $j\ne i$, since otherwise there would be a zero row (in $\mc{X}_{\Lambda\sm\{i\},\,\Lambda\sm\{n\}}(A)=A_{\Lambda\sm\{i\},\,\Lambda\sm\{n\}}$);
\item $j\in\Lambda\Rightarrow j+1\in\Lambda$ for all $j$ with $3\leq j\leq n-1$, since otherwise there would be a zero column.
\end{enumerate}

\noindent First assume that $i\ge3$ and that $|\Lambda|\le n-i+1$. Then it follows that $\Lambda=\{i,\ldots,n\}$ and that $\partial_{i n}\big(\det(\mc{X}_{\Lambda,\,\Lambda})\big)(A)=\pm1$.

Next assume that $i=2$. Then it follows that either $\Lambda=\{2,\ldots,n\}$ or $\Lambda=\{1,\ldots,n\}$. In the first case we have $\partial_{i n}\big(\det(\mc{X}_{\Lambda,\,\Lambda})\big)(A)=(-1)^{1+n-1}=(-1)^n$. In the second case we have $\partial_{i n}\big(\det(\mc{X}_{\Lambda,\,\Lambda})\big)(A)=(-1)^{2+n}=(-1)^n$.

Now assume that $i=1$. Then it follows that either $\Lambda=\{1,3,\ldots,n\}$ or $\Lambda=\{1,\ldots,n\}$. In the first case we have $\partial_{i n}\big(\det(\mc{X}_{\Lambda,\,\Lambda})\big)(A)=(-1)^{1+n-1}=(-1)^n$. In the second case we have $\partial_{i n}\big(\det(\mc{X}_{\Lambda,\,\Lambda})\big)(A)=(-1)^{1+n}.-1=(-1)^n$.

So for $i\in\{1,\ldots,n-1\}$ and $k\in\{1,\ldots,n\}$ we have:
$$(\partial_{i n}s_k)(A)=
\begin{cases}
\pm1 &\text{if $i\ge3$ and $i+k=n+1$,}\\
0 &\text{if $i\ge3$ and $i+k<n+1$,}\\
(-1)^n &\text{if $i\in\{1,2\}$ and $k\in\{n-1,n\}$,}\\
0 &\text{if $i\in\{1,2\}$ and $k<n-1$.}
\end{cases}
$$
It follows from the above equalities that in $M(A)$ the first $2$ columns are equal. So $d(A)=\det(M(A))=0$.

Let $\Lambda\subseteq\{1,\ldots,n\}$. Assume that $\partial_{1\,2}\big(\det(\mc{X}_{\Lambda,\,\Lambda})\big)$ is nonzero at $A$. Then $1,2\in\Lambda$ and the first row is zero. A contradiction. So $\partial_{1\,2}\big(\det(\mc{X}_{\Lambda,\,\Lambda})\big)$ is zero at $A$. Now assume that $\partial_{2\,1}\big(\det(\mc{X}_{\Lambda,\,\Lambda})\big)$ is nonzero at $A$. Then

\begin{enumerate}[$\bullet$]
\item $1,2\in\Lambda$;
\item $n\in\Lambda$, since otherwise the first row would be zero;
\item $j\in\Lambda\Rightarrow j-1\in\Lambda$ for all $j$ with $4\leq j\leq n$, since otherwise there would be a zero row.
\end{enumerate}
So $\Lambda=\{1,\ldots,n\}$ and $\partial_{i n}\big(\det(\mc{X}_{\Lambda,\,\Lambda})\big)(A)=\pm1$. Thus we have $(\partial_{1\,2}s_k)(A)=0$ for all $k\in\{1,\ldots,n\}$ and $(\partial_{2\,1}s_k)(A)=
\begin{cases}
\pm1&\text{if $k=n$,}\\
0 &\text{otherwise.}
\end{cases}$

Finally, we consider the function $d$. Let $i\in\{1,\ldots,n\}$, let $\Lambda\subseteq\{1,\ldots,n\}$ and assume that $\partial_{1\,2}\partial_{i n}\big(\det(\mc{X}_{\Lambda,\,\Lambda})\big)$ is nonzero at $A$. Then we have:
\begin{enumerate}[$\bullet$]
\item $1,2,i,n\in\Lambda$ and $i\ne 1$;
\item $i=2$, since otherwise the first row would be zero.
\item $j\in\Lambda\Rightarrow j-1\in\Lambda$ for all $j$ with $4\leq j\leq n$, since otherwise there would be a zero row.
\end{enumerate}
It follows that $i=2$, $\Lambda=\{1,\ldots,n\}$ and $\partial_{1\,2}\partial_{i n}\big(\det(\mc{X}_{\Lambda,\,\Lambda})\big)=\pm1$. So for $i,k\in\{1,\ldots,n\}$ we have:
$$(\partial_{1\,2}\partial_{i n}s_k)(A)=
\begin{cases}
\pm 1 &\text{if $(i,k)=(2,n)$,}\\
0 &\text{if $(i,k)\ne(2,n)$.}
\end{cases}
$$
We have
\begin{equation}\label{eq.d2}
d\,=\sum_{\pi\in {\mf S}_n} \sgn(\pi)\,\partial_{\pi(1),n}(s_1)\,\cdots\,\partial_{\pi(n),n}(s_n).
\end{equation}
So, by the above, $(\partial_{1\,2}d)(A)\,=$
$$\Big(\sum
\sgn(\pi)\partial_{\pi(1) n}(s_1)\partial_{\pi(2)n}(s_2)\cdots\partial_{\pi(n-1)n}(s_{n-1})
\partial_{1\,2}\partial_{2\,n}(s_n)\Big)(A),$$
where the sum is over all $\pi\in\mf{S}_n$ with $\pi(n)=2$. From what we know about the $\partial_{i n}s_k$ we deduce that the only permutation that survives in the above sum is given by $(\pi(1),\ldots,\pi(n))=(n,n-1,\ldots,3,1,2)$ and that $(\partial_{1\,2}d)(A)=\pm1$.

If we permute the rows of $M_\alpha(A)$ in the order given by $\Delta_1,\ldots,\Delta_{n-1},s_1,\ldots,s_n,d$ and take the columns in the order given by $\alpha$, then the resulting matrix is lower triangular with $\pm1$'s on the diagonal. So we can conclude that $d_\alpha(A)=\det(M_\alpha(A))=\pm1$.
\end{proof}

\noindent In the remainder of this subsection $K$ denotes an algebraically closed field.

It is well known that the algebra of regular functions $K[G]$ of a simply connected semi-simple algebraic group $G$ is a UFD (see \cite{Po} the corollary to Proposition~1), but the elementary proof below provides a way to show that $d'$ and the $\Delta_i'$ are irreducible elements of $K[\SL_n]$. I did not know how to use the fact that $K[\SL_n]$ is a UFD to simplify the proof that $\Delta_{n-1}'$ is irreducible.

Modifying the terminology of \cite{Eis} \S 16.6, we define the {\it Jacobian ideal} of an $m$-tuple of polynomials  $\varphi_1,\ldots,\varphi_m$ as the ideal generated by the $k\times k$ minors of the Jacobian matrix of $\varphi_1,\ldots,\varphi_m$, where $k$ is the height of the ideal generated by the $\varphi_i$.

\begin{lemgl}\label{lem.sln}
$K[\SL_n]$ is a unique factorisation domain and $\Delta'_{n-1}$ is an irreducible element of $K[\SL_n]$.
\end{lemgl}

\begin{proof}
From the Laplace expansion for $\det$ with respect to the last row or the last column it is clear that we can eliminate $\xi_{nn}$ using the relation $\det=1$, if we make $\Delta'_{n-1}$ invertible. So we have an isomorphism of $K[\SL_n][\Delta_{n-1}'^{\ -1}]$ with the localised polynomial algebra $K[(\xi_{ij})_{(i,j)\ne(n,n)}][\Delta_{n-1}^{-1}]$. Since the latter algebra is a UFD, it suffices, by Nagata's lemma, to prove that $\Delta'_{n-1}$ is prime in $K[\SL_n]$, i.e. that $(\Delta_{n-1},\det-1)$ generates a prime ideal in $K[\Mat_n]$. First we show that the closed subvariety $\mc V$ of $\Mat_n$ defined by this ideal is irreducible.

Let $\mc{X}$ the matrix introduced above and let $\alpha_1,\ldots,\alpha_{n-2}$ be variables. For a matrix $A$ denote by $A^{(i,j)}$ the matrix which is obtained from $A$ by deleting the $i$-th row and the $j$-th column. Let ${\mc X}_\alpha$ be the $n\times n$ matrix which is obtained by replacing in $\mc X$ the $(n-1)$-th column of ${\mc X}^{(n,n)}$ by the linear combination $\sum_{j=1}^{n-2}\alpha_j(\xi_{ij})_i$ of the first $n-2$ columns of ${\mc X}^{(n,n)}$. Then $\det({\mc X}_\alpha^{(n,j)})=\pm \alpha_j\det({\mc X}^{(n,n-1)})$ for all $j\in\{1,\ldots,n-2\}$ and $\det({\mc X}_\alpha^{(n,n)})=0$. So, by the Laplace expansion rule
\begin{align*}
\det({\mc X}_\alpha)-1&=\sum_{j=1}^{n-1}\pm\xi_{nj}\det({\mc X_\alpha}^{(n,j)})-1\\
&=\pm\xi_{n\,n-1}\det({\mc X}^{(n,n-1)})+\sum_{j=1}^{n-2}\pm\alpha_j\xi_{nj}\det({\mc X}^{(n,n-1)})-1
\end{align*}
Let $K[{\mc X}_\alpha]$ be the polynomial ring in the variables that occur in ${\mc X}_\alpha$. If we consider $\det({\mc X}_\alpha)-1$ as a polynomial in the variable $\xi_{n,n-1}$, then it is linear and its leading coefficient is $\pm\det({\mc X}^{(n,n-1)})$ which is irreducible and does not divide the constant term $\sum_{j=1}^{n-2}\pm\alpha_j\xi_{nj}\det({\mc X}^{(n,n-1)})-1$. So $\det({\mc X}_\alpha)-1$ is irreducible in $K[{\mc X}_\alpha]$ and it defines an irreducible closed subvariety ${\mc V}_1$ of an $n^2-1$ dimensional affine space with coordinate functionals $\xi_{ij}$, $j\ne n-1$, $\xi_{n\,n-1}$, $\alpha_1,\ldots,\alpha_{n-2}$.

Let $H$ be the algebraic group of $n\times n$ matrices $(a_{ij})$ of determinant $1$ with $a_{nn}=1$ and $a_{in}=a_{ni}=0$ for all $i\in\{1,\ldots,n-1\}$. Then $H\cong\SL_{n-1}$ and for every $A\in{\mc V}$ there exists an $S\in H$ such that in $(AS)^{(n,n)}$ the last column is a linear combination of the others. So the morphism of varieties $(u,S)\mapsto {\mc X}_\alpha(u)S:{\mc V}_1\times H\to\Mat_n$ has image $\mc V$. Now the irreducibility of $\mc V$ follows from the irreducibility of ${\mc V}_1\times H$.

It remains to show that $(\det-1,\Delta_{n-1})$ is a radical ideal of $K[\Mat_n]$, i.e. that $K[\Mat_n]/(\det-1,\Delta_{n-1})$ is reduced. We know that $\Delta_{n-1}\ne0$ on the irreducible variety $\SL_n$, so $\dim({\mc V})=n^2-2$ and $K[\Mat_n]/(\det-1,\Delta_{n-1})$ is Cohen-Macaulay (see \cite{Eis} Proposition 18.13). By Theorem 18.15 in \cite{Eis} it suffices to show that the closed subvariety of $\mc V$ defined by the Jacobian ideal of $\det-1,\Delta_{n-1}$ is of codimension $\ge1$. Since $\mc V$ is irreducible this follows from Lemma~\ref{lem.jac}(i).
\end{proof}

\begin{lemgl}\label{lem.slnd}
\begin{enumerate}[(i)]
\item $d$ is an irreducible element of $K[\Mat_n]$.
\item The invertible elements of $K[\SL_n]$ are the nonzero scalars.
\item $d',\Delta'_1,\ldots,\Delta'_{n-1}$ is are mutually inequivalent irreducible elements of $K[\SL_n]$.
\end{enumerate}
\end{lemgl}

\begin{proof}
(i).\ The proof of this is completely analogous to that of Proposition 3 in \cite{PrT}. One now has to work with the maximal parabolic subgroup $P$ of $\GL_n$ that consists of the invertible matrices $(a_{ij})$ with $a_{n\,i}=0$ for all $i<n$. The element $d$ is then a semi-invariant of $P$ with the weight $\det\cdot\xi_{nn}^{-n}$ (the restriction of this weight to the maximal torus of diagonal matrices is $n\varpi_{n-1}$).\\
(ii) and (iii).\ Consider the isomorphism $K[\SL_n][\Delta_{n-1}'^{\ -1}]\cong K[(\xi_{ij})_{(i,j)\ne(n,n)}][\Delta_{n-1}^{-1}]$ from the proof of the above lemma. It maps $d',\Delta'_1,\ldots,\Delta'_{n-1}$ to respectively $d,\Delta_1,\ldots,\Delta_{n-1}$, since these polynomials do not contain the variable $\xi_{nn}$. The invertible elements of $K[(\xi_{ij})_{(i,j)\ne(n,n)}][\Delta_{n-1}^{-1}]$ are the elements $\alpha\Delta_{n-1}^k$, $\alpha\in K\sm\{0\}$, $k\in{\mb Z}$, since $\Delta_{n-1}$ is irreducible in $K[(\xi_{ij})_{(i,j)\ne(n,n)}]$. So the invertible elements of $K[\SL_n][\Delta_{n-1}'^{\ -1}]$ are the elements $\alpha\Delta_{n-1}'^{\ k}$, $\alpha\in K\sm\{0\}$, $k\in{\mb Z}$. By Lemma~\ref{lem.sln} $\Delta_{n-1}'$ is irreducible in $K[\SL_n]$, so the invertible elements of $K[\SL_n]$ are the nonzero scalars. Since $d$ and the $\Delta_i$ are not scalar multiples of each other, all that remains is to show that $d'$ and $\Delta'_1,\ldots,\Delta'_{n-2}$ are irreducible. We only do this for $d'$, the argument for the $\Delta_i'$ is completely similar.
Since $d$ is prime in $K[(\xi_{ij})_{(i,j)\ne(n,n)}]$ and $d$ does not divide $\Delta_{n-1}$, it follows that $d$ is prime in $K[(\xi_{ij})_{(i,j)\ne(n,n)}][\Delta_{n-1}^{-1}]$ and therefore that $d'$ is prime in $K[\SL_n][\Delta_{n-1}'^{\ -1}]$. To show that $d'$ is prime in $K[\SL_n]$ it suffices to show that for every $f\in K[\SL_n]$, $\Delta'_{n-1}f\in(d')$ implies $f\in(d')$. So assume that $\Delta'_{n-1}f=gd'$ (*) for some $f,g\in K[\SL_n]$. If we take ${\bf a}\in K^n$ such that $a_n=(-1)^{n-1}$, then we have $x_{\bf a}\in\SL_n$, $d'(x_{\bf a})=1$ and $\Delta'_{n-1}(x_{\bf a})=0$. So $\Delta'_{n-1}$ does not divide $d'$. But then, by Lemma~\ref{lem.sln}, $\Delta'_{n-1}$ divides $g$. Cancelling a factor $\Delta'_{n-1}$ on both sides of (*), we obtain that $f\in(d')$.
\end{proof}

\subsection{Generators and relations and a $\mb Z$-form for $\tilde{Z_0}[z_{\varpi_1},\ldots,z_{\varpi_{n-1}}]Z_1$}\label{ss.presentation}
\ \\
For the basics about monomial orderings and Gr\"obner bases I refer to \cite{CLO}.

\begin{lemgl}\label{lem.leadmondetd}
If we give the monomials in the variables $\xi_{ij}$ the lexicographic monomial ordering for which $\xi_{n\,n}>\xi_{n\,n-1}\cdots>\xi_{n1}>\xi_{n-1\,n}>\cdots> \xi_{n-1\,1}>\cdots>\xi_{11}$, then $\det$ has leading term $\pm\xi_{n\,n}\cdots\xi_{2\,2}\xi_{1\,1}$ and $d$ has leading term $\pm\xi_{n\,n-1}^{n-1}\cdots\xi_{3\,2}^2\xi_{2\,1}$.
\end{lemgl}
\newcommand{\LT}{{\rm LT}}
\begin{proof}
\newcommand{\pol}{\partial_{i\,n}\big(\det(\mc{X}_{\Lambda,\,\Lambda})\big)}
I leave the proof of the first assertion to the reader. For the second assertion we use the notation and the formulas of Subsection~\ref{ss.slnd}. The leading term of a nonzero polynomial $f$ is denoted by $\LT(f)$. Let $i\in\{1,\ldots,n\}$ and $\Lambda\subseteq\{1,\ldots,n\}$ with $|\Lambda|=k\ge2$ and assume that $\pol\ne0$. Then $i,n\in\Lambda$. Now we use the fact that no monomial in $\pol$ contains a variable with row index equal to $i$ or with column index equal to $n$ or a product of two variables which have the same row or column index.

First assume that $i>n-k+1$. Then $$\LT(\pol)\le\pm\xi_{n\,n-1}\cdots\xi_{i+1\,i}\xi_{i-1\,i-1}\cdots\xi_{n-k+1\,n-k+1}$$ with equality if and only if $\Lambda=\{n,n-1,\ldots,n-k+1\}$. Now assume that $i=n-k+1$. Then $$\LT(\pol)\le\pm\xi_{n\,n-1}\cdots\xi_{n-k+2\,n-k+1}$$ with equality if and only if $\Lambda=\{n,n-1,\ldots,n-k+1\}$. Finally assume that $i<n-k+1$. Then
$$\LT(\pol)\le\pm\xi_{n\,n-1}\cdots\xi_{n-k+3\,n-k+2}\xi_{n-k+2\,i}$$
with equality if and only if $\Lambda=\{n,n-1,\ldots,n-k+2,i\}$.

So for $i,k\in\{1,\ldots,n\}$ with $k\ge2$ we have:
$$\LT(\partial_{i n}s_k)=
\begin{cases}
\pm\xi_{n\,n-1}\cdots\xi_{i+1\,i}\xi_{i-1\,i-1}\cdots\xi_{n-k+1\,n-k+1} &\text{if $i+k>n+1$,}\\
\pm\xi_{n\,n-1}\cdots\xi_{n-k+2\,n-k+1} &\text{if $i+k=n+1$,}\\
\pm\xi_{n\,n-1}\cdots\xi_{n-k+3\,n-k+2}\xi_{n-k+2\,i} &\text{if $i+k<n+1$.}
\end{cases}
$$
In particular $\LT(\partial_{i n}s_k)\le\pm\xi_{n\,n-1}\cdots\xi_{n-k+1\,n-k+1}$ with equality if and only if $i+k=n+1$. But then, by equation~\eqref{eq.d2}, $\LT(d)=\LT(\partial_{n\,n}s_1)\LT(\partial_{n-1\,n}s_2)\cdots\LT(\partial_{1\,n}s_n)$\\
$=\pm\xi_{n\,n-1}^{n-1}\cdots\xi_{3\,2}^2\xi_{2\,1}$
\end{proof}

Recall that the degree reverse lexicographical ordering on the monomials $u^\alpha=u_1^{\alpha_1}\cdots u_k^{\alpha_k}$ in the variables $u_1,\ldots,u_k$ is defined as follows: $u^\alpha > u^\beta$ if $\deg(u^\alpha)>\deg(u^\beta)$ or $\deg(u^\alpha)=\deg(u^\beta)$ and $\alpha_i<\beta_i$ for the last index $i$ with $\alpha_i\ne\beta_i$.

\begin{lemgl}\label{lem.leadmonfi}
Let $f_i\in\mb{Z}[u_1,\ldots,u_{n-1}]$ be the polynomial such that $\sym(l\varpi_i)=$\\ $f_i(\sym(\varpi_1),\ldots,\sym(\varpi_{n-1}))$. If we give the monomials in the $u_i$ the degree reverse lexicographic monomial ordering for which $u_1>\cdots>u_{n-1}$, then $f_i$ has leading term $u_i^l$. Furthermore, the monomials that appear in $f_i-u_i^l$ are of total degree $\le l$ and have exponents $<l$.\ \footnotemark
\end{lemgl}

\footnotetext{So our $f_i$ are related to the polynomials $P_i=x_i^l-\sum_\mu d_{i\mu}x_\mu$ from the proof of Proposition 6.4 in \cite{DCKP} as follows: $P_i=f_i(x_1,\ldots,x_{n-1})-\sym(l\varpi_i)$. In particular $d_{i0}=\sym(l\varpi_i)$ and $d_{i\mu}\in{\mb Z}$ for all $\mu\in P\sm\{0\}$ (we are, of course, in the situation that $\g=\spl_n$).}

\begin{proof}
Let $\sigma_i$ be the $i$-th elementary symmetric function in the variables $x_1,\ldots,x_n$ and let $\lambda_i\in P=X(T)$ be the character $A\mapsto A_{ii}$ of $T$. Then $\sym(\varpi_i)=\sigma_i(e(\lambda_1),\ldots,e(\lambda_n))$ for $i\in\{1,\ldots,n-1\}$. So the $f_i$ can be found as follows. For $i\in\{1,\ldots,n-1\}$, determine $F_i\in\mb{Z}[u_1,\ldots,u_n]$ such that $\sigma_i(x_1^l,\ldots,x_n^l)=F_i(\sigma_1,\ldots,\sigma_n)$. Then $f_i=F_i(u_1,\ldots,u_{n-1},1)$. It now suffices to show that for $i\in\{1,\ldots,n-1\}$, $F_i-u_i^l$ is a $\mb Z$-linear combination of monomials in the $u_j$ that have exponents $<l$, are of total degree $\le l$ and that contain some $u_j$ with $j>i$ (the monomials that contain $u_n$ will become of total degree $<l$ when $u_n$ is replaced by 1).

Fix $i\in\{1,\ldots,n-1\}$. Consider the following properties of a monomial in the $x_j$:
\begin{enumerate}[(x1)]
\item the monomial contains at least $i+1$ variables.
\item the exponents are $\le l$.
\item the number of exponents equal to $l$ is $\le i$.
\end{enumerate}
and the following properties of a monomial in the $u_j$:
\begin{enumerate}[(u1)]
\item the monomial contains a variable $u_j$ for some $j>i$.
\item the total degree is $\le l$.
\item the exponents are $<l$.
\end{enumerate}
Let $h$ be a symmetric polynomial in the $x_i$ and let $H$ be the polynomial in the $u_i$ such that $h=H(\sigma_1,\ldots,\sigma_n)$. Give the monomials in the $x_i$ the lexicographic monomial ordering for which $x_1>\cdots>x_n$. We will show by induction on the leading monomial of $h$ that if each monomial that appears in $h$ has property (x1) resp. property (x2) resp. properties (x1), (x2) and (x3), then each monomial that appears in $H$ has property (u1) resp. property (u2) resp. properties (u1), (u2) and (u3). Let $x^\alpha:=x_1^{\alpha_1}\cdots x_n^{\alpha_n}$ be the leading monomial of $h$. Then $\alpha_1\ge\alpha_2\cdots\ge\alpha_n$. Put $\beta=(\alpha_1-\alpha_2,\ldots,\alpha_{n-1}-\alpha_n,\alpha_n)$. Let $k$ be the last index for which $\alpha_k\ne 0$. Then $\beta=(\alpha_1-\alpha_2,\ldots,\alpha_{k-1}-\alpha_k,\alpha_k,0,\ldots,0)$. If $x^\alpha$ has property (x1), then $k\ge i+1$, $u^\beta$ has property (u1) and the monomials that appear in $\sigma^\beta$ have property (x1), since $\sigma_k$ appears in $\sigma^\beta$. 
If $x^\alpha$ has property (x2), then $\alpha_1\le l$, $u^\beta$ is of total degree $\alpha_1\le l$ and the monomials that appear in $\sigma^\beta$ have exponents $\le\beta_1+\cdots+\beta_k=\alpha_1\le l$. Now assume that $x^\alpha$ has properties (x1), (x2) and (x3). For $j<k$ we have $\beta_j=\alpha_j-\alpha_{j+1}<l$, since $\alpha_{j+1}\ne0$. So we have to show that $\beta_k=\alpha_k<l$. If $\alpha_k$ were equal to $l$, then we would have $\alpha_1=\cdots=\alpha_k=l$, by (x2). This contradicts (x3), since we have $k\ge i+1$ by (x1). Finally we show that the monomials that appear in $\sigma^\beta$ have property (x3). If $\alpha_1<l$, then all these monomials have exponents $<l$. So assume $\alpha_1=l$. Let $j$ be the smallest index for which $\beta_j\ne 0$. Then the number of exponents equal to $l$ in a monomial that appears in $\sigma^\beta$ is $\le j$. On the other hand $\alpha_1=\cdots=\alpha_j=l$. So we must have $j\le i$, since $x^\alpha$ has property (x3).

Now we can apply the induction hypothesis to $h-c\sigma^\beta$, where $c$ is the leading coefficient of $h$.

The assertion about $F_i-u_i^l$ now follows, because the monomials that appear in $\sigma_i(x_1^l,\ldots,x_n^l)-\sigma_i^l$ have the properties (x1), (x2) and (x3).
\end{proof}

From now on we denote $z_{\varpi_i}$ by $z_i$.\footnote{In \cite{DCKP} and \cite{DCP} $z_{\alpha_i}$ is denoted by $z_i$.} Let $\mb{Z}[\SL_n]$ be the $\mb Z$-subalgebra of $\mb{C}[\SL_n]$ generated by the $\xi_{ij}'$ and $A$ be the $\mb Z$-subalgebra of $Z$ generated by the $\tilde{\xi}_{ij}$. So $A=\pi^{co}(\mb{Z}[\SL_n])$. Let $B$ be the $\mb Z$-subalgebra generated by the elements $\tilde{\xi}_{ij}$, $u_1,\ldots,u_{n-1}$ and $z_1,\ldots,z_{n-1}$. For a commutative ring $R$ we put $A(R)=R\ot_{\mb Z}A$ and $B(R)=R\ot_{\mb Z}B$. Clearly we can identify $A({\mb C})$ with $\tilde{Z_0}$. In the proposition below "natural homomorphism" means a homomorphism that maps $\xi_{ij}$ to $\tilde{\xi}_{ij}$ and, if this applies, the variables $u_i$ and $z_i$ to the equally named elements of $Z$. The polynomials $f_i$ below are the ones defined in Lemma~\ref{lem.leadmonfi}.

\begin{propgl}\label{prop.presentation}
The following holds:
\begin{enumerate}[(i)]
\item The kernel of the natural homomorphism from the polynomial algebra\\ $\mb{Z}[(\xi_{ij})_{ij},u_1,\ldots,u_{n-1},z_1,\ldots,z_{n-1}]$ to $B$ is generated by the elements\\ $\det-1,f_1-s_1,\ldots,f_{n-1}-s_{n-1},z_1^2-\Delta_1,\ldots,z_{n-1}^2-\Delta_{n-1}$.
\item The homomorphism $B(\mb{C})\to Z$, given by the universal property of ring transfer, is injective.
\item $A$ is a free $\mb{Z}$-module and $B$ is a free $A$-module with the monomials\\
$u_1^{k_1}\cdots u_{n-1}^{k_{n-1}}z_1^{m_1}\cdots z_{n-1}^{m_{n-1}}$, $0\leq k_i<l$, $0\leq m_i<2$ as a basis.
\item $A[z_1,\ldots,z_{n-1}]\cap Z_1=A\cap Z_1={\mb Z}[\tilde{s}_1,\ldots,\tilde{s}_{n-1}]$ and $B\cap Z_1$ is a free $A\cap Z_1$-module with the monomials $u_1^{k_1}\cdots u_{n-1}^{k_{n-1}}$, $0\leq k_i<l$ as a basis.
\end{enumerate}
\end{propgl}

\begin{proof}
Let $Z_0'$ be the $\mb C$-subalgebra of $Z$ generated by the $\tilde{\xi}_{ij}$ and $z_1,\ldots,z_{n-1}$. As we have seen in Subsection~\ref{ss.Z0andZ1}, the $z_i$ satisfy the relations $z_i^2=\tilde{\Delta_i}$. The $\tilde{\Delta_i}$ are part of a generating transcendence basis of the field of fractions ${\rm Fr}(\tilde{Z_0})$ of $\tilde{Z_0}$ by arguments very similar to those at the end of the proof of Theorem~\ref{thm.qrat}. This shows that the monomials $z_1^{m_1}\cdots z_{n-1}^{m_{n-1}}$, $0\leq m_i<2$, form a basis of ${\rm Fr}(Z_0')$ over ${\rm Fr}(\tilde{Z_0})$ and of $Z_0'$ over $\tilde{Z_0}$. It follows that the kernel of the natural homomorphism from the polynomial algebra $\mb{C}[(\xi_{ij})_{ij},z_1,\ldots,z_{n-1}]$ to $Z_0'$ is generated by the elements $\det-1,z_1^2-\Delta_1,\ldots,z_{n-1}^2-\Delta_{n-1}$. So we have generators and relations for $Z_0'$. By the construction from Subsection~\ref{ss.Z0andZ1} we then obtain that the kernel $I$ of the natural homomorphism from the polynomial algebra $\mb{C}[(\xi_{ij})_{ij},u_1,\ldots,u_{n-1},z_1,\ldots,z_{n-1}]$ to $Z_0'Z_1$ is generated by the elements $\det-1,f_1-s_1,\ldots,f_{n-1}-s_{n-1},z_1^2-\Delta_1,\ldots,z_{n-1}^2-\Delta_{n-1}$.

Now we give the monomials in the variables $(\xi_{ij})_{ij},u_1,\ldots,u_{n-1},z_1,\ldots,z_{n-1}$ a monomial ordering which is the lexicographical product of an arbitrary monomial ordering on the monomials in the $z_i$, the monomial ordering of Lemma~\ref{lem.leadmonfi} on the monomials in the $u_i$ and the monomial ordering of Lemma~\ref{lem.leadmondetd} on the $\xi_{ij}$.\ \footnote{so the $z_i$ are greater than the $u_i$ which are greater than the $\xi_{ij}$} Then the ideal generators mentioned above have leading monomials $\xi_{n\,n}\cdots\xi_{2\,2}\xi_{1\,1},$
$u_1^l,\ldots,u_{n-1}^l, z_1^2,\ldots,z_{n-1}^2$ and the leading coefficients are all $\pm1$. Since the leading monomials have gcd $1$, the ideal generators form a Gr\"obner basis; see \cite{CLO} Ch.~2 \S~9 Theorem~3 and Proposition~4, for example. Since the leading coefficients are all $\pm1$, it follows from the division with remainder algorithm that the ideal of $\mb{Z}[(\xi_{ij})_{ij},u_1,\ldots,u_{n-1},z_1,\ldots,z_{n-1}]$ generated by these elements consists of the polynomials in $I$ that have integral coefficients and that it has the $\mb Z$-span of the monomials that are not divisible by any of the above leading monomials as a direct complement. This proves (i) and (ii).\\
(iii).\ The canonical images of the above monomials form a $\mb Z$-basis of $B$. These monomials are the products of the monomials in the $\xi_{ij}$ that are not divisible by $\xi_{n\,n}\cdots\xi_{2\,2}\xi_{1\,1}$ and the restricted monomials mentioned in the assertion. The canonical images of the monomials in the $\xi_{ij}$ that are not divisible by $\xi_{n\,n}\cdots\xi_{2\,2}\xi_{1\,1}$ form a $\mb Z$-basis of $A$.\\
(iv).\ As we have seen, the monomials with exponents $<2$ in the $z_i$ form a basis of the $\tilde{Z_0}$-module $Z_0'$. So $A[z_1,\ldots,z_{n-1}]\cap \tilde{Z_0}=A$. Therefore, by Theorem~\ref{thm.coverHCcentre}(ii), $A[z_1,\ldots,z_{n-1}]\cap Z_1=A\cap Z_1=\pi^{co}({\mb Z}[\SL_n]^{\SL_n})$. Now $({\mb Z}P)^W={\mb Z}[\sym(\varpi_1),\ldots,\sym(\varpi_{n-1})]$ (see \cite{Bou1} no. VI.3.4, Thm. 1.) and the $s_i'$ are in ${\mb Z}[\SL_n]$, so ${\mb Z}[\SL_n]^{\SL_n}={\mb Z}[s_1',\ldots,s_{n-1}']$ by the restriction theorem for ${\mb C}[\SL_n]$. This proves the first assertion. From the proof of Theorem~\ref{thm.Z0andZ1} we know that the given monomials form a basis of $Z_1$ over $Z_0\cap Z_1$ and a basis of $Z$ over $Z_0$. So an element of $Z$ is in $Z_1$ if and only if its coefficients with respect to this basis are in $Z_0\cap Z_1$. The second assertion now follows from (iii).
\end{proof}

By (ii) of the above proposition we may identify $B({\mb C})$ with $\tilde{Z_0}[z_1,\ldots,z_{n-1}]Z_1$ and $B({\mb C})[\tilde{\Delta}^{-1}_1,\ldots,\tilde{\Delta}^{-1}_{n-1}]$ with $Z$.

Put $\ov{Z}=Z/(\tilde{d})$. For the proof of Theorem~\ref{thm.quf} we need a version for $\ov{Z}$ of Proposition~\ref{prop.presentation}. First we introduce some more notation. For $u\in Z$ we denote the canonical image of $u$ in $\ov{Z}$ by $\ov{u}$. For $f\in\mb{C}[\Mat_n]$ we write $\ov{f}$ instead of $\ov{\tilde{f}}$. Let $\ov{A}$ be the $\mb Z$-subalgebra of $\ov{Z}$ generated by the $\ov{\xi}_{ij}$ and let $\ov{B}$ be the $\mb Z$-subalgebra generated by the elements $\ov{\xi}_{ij}$, $\ov{u}_1,\ldots,\ov{u}_{n-1}$ and $\ov{z}_1,\ldots,\ov{z}_{n-1}$. For a commutative ring $R$ we put $\ov{A}(R)=R\ot_{\mb Z}\ov{A}$ and $\ov{B}(R)=R\ot_{\mb Z}\ov{B}$.

\renewcommand{\thepropnn}{$\ov{\text{\ref{prop.presentation}}}$}
\begin{propnn}\label{prop.presentation(d)}
The following holds:
\begin{enumerate}[(i)]
\item The kernel of the natural homomorphism from the polynomial algebra\\ $\mb{Z}[(\xi_{ij})_{ij},u_1,\ldots,u_{n-1},z_1,\ldots,z_{n-1}]$ to $\ov{B}$ is generated by the elements\\ $\det-1,d,f_1-s_1,\ldots,f_{n-1}-s_{n-1},z_1^2-\Delta_1,\ldots,z_{n-1}^2-\Delta_{n-1}$.
\item The kernel of the natural homomorphism $\mb{Z}[\Mat_n]\to\ov{A}$ is $(\det-1,d)$.
\item The homomorphism $\ov{B}(\mb{C})\to \ov{Z}$, given by the universal property of ring transfer, is injective.
\item $\ov{A}$ is a free $\mb{Z}$-module and $\ov{B}$ is a free $\ov{A}$-module with the monomials\\ $\ov{u}_1^{k_1}\cdots \ov{u}_{n-1}^{k_{n-1}}\ov{z}_1^{m_1}\cdots \ov{z}_{n-1}^{m_{n-1}}$, $0\leq k_i<l$, $0\leq m_i<2$ as a basis.
\item The $\ov{A}$-span of the monomials $\ov{u}_1^{k_1}\cdots \ov{u}_{n-1}^{k_{n-1}}$, $0\leq k_i<l$, is closed under multiplication.
\end{enumerate}
\end{propnn}
\renewcommand{\thepropnn}{\!\!}
\begin{proof}
From Lemma~\ref{lem.slnd}(iii) we deduce that $\big(A({\mb C})[\tilde{\Delta}_1^{-1},\ldots,\tilde{\Delta}_{n-1}^{-1}]\tilde{d}\big)\cap A({\mb C})=A({\mb C})\tilde{d}$. From this it follows, using the $A({\mb C})$-basis of $B({\mb C})$, that $(Z\tilde{d})\cap B({\mb C})$, which is the kernel of the natural homomorphism $B({\mb C})\to\ov{Z}$, equals $B({\mb C})\tilde{d}$. From (i) and (ii) of Proposition~\ref{prop.presentation} or from its proof it now follows that the kernel of the natural homomorphism from the polynomial algebra\\ $\mb{C}[(\xi_{ij})_{ij},u_1,\ldots,u_{n-1},z_1,\ldots,z_{n-1}]$ to $\ov{Z}$ is generated by the elements $\det-1,d,f_1-s_1,\ldots,f_{n-1}-s_{n-1},z_1^2-\Delta_1,\ldots,z_{n-1}^2-\Delta_{n-1}$.

Again using the $A({\mb C})$-basis of $B({\mb C})$ we obtain that $(B({\mb C})\tilde{d})\cap A({\mb C})=A({\mb C})\tilde{d}$. From this it follows that the kernel of the natural homomorphism $\mb{C}[\Mat_n]\to\ov{Z}$ is generated by $\det-1$ and $d$.

By Lemma~\ref{lem.leadmondetd} we have $\LT(d)=\pm\xi_{n\,n-1}^{n-1}\cdots\xi_{3\,2}^2\xi_{2\,1}$ which has gcd $1$ with the leading monomials of the other ideal generators, so the ideal generators mentioned above form a Gr\"obner basis over $\mb Z$. Now (i)-(iv) follow as in the proof of Proposition~\ref{prop.presentation}.\\
(v).\ This follows from the fact that the remainder modulo the Gr\"obner basis of a polynomial in $\mb{Z}[(\xi_{ij})_{ij},u_1,\ldots,u_{n-1}]$ is again in $\mb{Z}[(\xi_{ij})_{ij},u_1,\ldots,u_{n-1}]$.
\end{proof}

By (ii) and (iii) of the above proposition $\ov{A}$ and $\ov{B}({\mb C})[\ov{\Delta}^{\ -1}_1,\ldots,\ov{\Delta}^{\ -1}_{n-1}]$ can be identified with respectively ${\mb Z}[\Mat_n]/(\det-1,d)$ and $\ov{Z}$. From (iv) it follows that, for any commutative ring $R$, $\ov{A}(R)$ embeds in $\ov{B}(R)$.

\subsection{The theorem}

\begin{lemgl}\label{lem.algnum}
Let $A$ be an associative algebra with 1 over a field $F$ and let $L$ be an extension of $F$. Assume that for every finite extension $F'$ of $F$, $F'\ot_FA$ has no zero divisors. Then the same holds for $L\ot_FA$.
\end{lemgl}

\begin{proof}
Assume that there exist $a,b\in L\ot_FA\sm\{0\}$ with $ab=0$. Let $(e_i)_{i\in I}$ be an $F$-basis of $A$ and let $c_{ij}^k\in F$ be the structure constants. Write $a=\sum_{i\in I}\alpha_ie_i$ and $b=\sum_{i\in I}\beta_ie_i$. Let $I_a$ resp. $I_b$ be the set of indices $i$ such that $\alpha_i\ne0$ resp. $\beta_i\ne0$ and let $J$ be the set of indices $k$ such that $c_{ij}^k\ne0$ for some $(i,j)\in I_a\times I_b$. Then $I_a$ and $I_b$ are nonempty and $I_a$, $I_b$ and $J$ are finite. Take $i_a\in I_a$ and $i_b\in I_b$. Since $ab=0$, the following equations over $F$ in the variables $x_i$, $i\in I_a$, $y_i$, $i\in I_b$, $u$ and $v$ have a solution over $L$:
\begin{align*}
&\sum_{i\in I_a, j\in I_b}c_{ij}^kx_iy_j=0\text{ for all }k\in J,\\
&x_{i_a}u=1, y_{i_b}v=1.
\end{align*}

But then they also have a solution over a finite extension $F'$ of $F$ by Hilbert's Nullstellensatz. This solution gives us nonzero elements $a',b'\in F'\ot_FA$ with $a'b'=0$.
\end{proof}

\begin{lemgl}\label{lem.modp}
Let $R$ be the valuation ring of a nontrivial discrete valuation of a field $F$ and let $K$ be its residue class field. Let $A$ be an associative algebra with 1 over $R$ which is free as an $R$-module and let $L$ be an extension of $F$. Assume that for every finite extension $K'$ of $K$, $K'\ot_RA$ has no zero divisors. Then the same holds for $L\ot_RA$.
\end{lemgl}

\begin{proof}
Assume that there exist $a,b\in L\ot_RA\sm\{0\}$ with $ab=0$. By the above lemma we may assume that $a,b\in F'\ot_RA\sm\{0\}$ for some finite extension $F'$ of $F$. Let $(e_i)_{i\in I}$ be an $R$-basis of $A$. Let $\nu$ be an extension to $F'$ of the given valuation of $F$, let $R'$ be the valuation ring of $\nu$, let $K'$ be the residue class field and let $\delta\in R'$ be a uniformiser for $\nu$. Note that $R'$ is a local ring and a principal ideal domain (and therefore a UFD) and that $K'$ is a finite extension of $K$ (see e.g. \cite{Cohn} Chapter 8 Theorem 5.1). By multiplying $a$ and $b$ by suitable integral powers of $\delta$ we may assume that their coefficients with respect to the basis $(e_i)_{i\in I}$ are in $R'$ and not all divisible by $\delta$ (in $R'$). By passing to the residue class field $K$ we then obtain nonzero $a',b'\in K'\ot_{R'}(R'\ot_RA)=K'\ot_RA$ with $a'b'=0$.
\end{proof}

\begin{remnn}
The above lemmas also hold if we replace "zero divisors" by "nonzero nilpotent elements".
\end{remnn}

For $t\in\{0,\ldots,n-1\}$ let $\ov{B}_t$ be the $\mb Z$-subalgebra generated by the elements $\ov{\xi}_{ij}$, $\ov{u}_1,\ldots,\ov{u}_{n-1}$ and $\ov{z}_1,\ldots,\ov{z}_t$. So $\ov{B}_{n-1}=\ov{B}$. For a commutative ring $R$ we put $\ov{B}_t(R)=R\ot_{\mb Z}\ov{B}_t$. From (iv) and (v) of Proposition~\ref{prop.presentation(d)} we deduce that the monomials $\ov{u}_1^{k_1}\cdots \ov{u}_{n-1}^{k_{n-1}}\ov{z}_1^{m_1}\cdots \ov{z}_t^{m_t}$, $0\leq k_i<l$, $0\leq m_i<2$ form a basis of $\ov{B}_t$ over $\ov{A}$. So for any commutative ring $R$ we have bases for $\ov{B}_t(R)$ over $\ov{A}(R)$ and over $R$. Note that $\ov{B}_t(R)$ embeds in $\ov{B}(R)$, since the $\mb Z$-basis of $\ov{B}_t$ is part of the $\mb Z$-basis of $\ov{B}$.

\begin{thmgl}\label{thm.quf}
If $l$ is a power of an odd prime $p$, then $Z$ is a unique factorisation domain.
\end{thmgl}
\begin{proof}

We have seen in Subsection~\ref{ss.nis2} that for $n=2$ it holds without any extra assumptions on $l$, so assume that $n\ge3$. For the elimination of variables in the proof of Theorem~\ref{thm.qrat} we only needed the invertibility of $\tilde{d}$, so $Z[\tilde{d}^{-1}]$ is isomorphic to a localisation of a polynomial algebra and therefore a UFD. So, by Nagata's lemma, it suffices to prove that $\tilde{d}$ is a prime element of $Z$, i.e. that $\ov{Z}=Z/(\tilde{d})$ is an integral domain. We do this in 5 steps.

\noindent 1.\ $\ov{B}(K)$ is reduced for any field $K$.

We may assume that $K$ is algebraically closed. Since $\ov{B}(K)$ is a finite $\ov{A}(K)$-module it follows that $\ov{B}(K)$ is integral over $\ov{A}(K)\cong K[\Mat_n]/(\det-1,d)$. So it its Krull dimension is $n^2-2$. By Proposition~\ref{prop.presentation(d)}, $\ov{B}(K)$ is isomorphic to the quotient of a polynomial ring over $K$ in $n^2+2(n-1)$ variables by an ideal $I$ which is generated by $2n$ elements. So $\ov{B}(K)$ is Cohen-Macaulay (see \cite{Eis} Proposition 18.13). Let $\mc V$ be the closed subvariety of $n^2+2(n-1)$-dimensional affine space defined by $I$. By Theorem 18.15 in \cite{Eis} it suffices to show that the closed subvariety of $\mc V$ defined by the Jacobian ideal of $\det-1,d, f_1-s_1,\ldots,f_{n-1}-s_{n-1}, z_1^2-\Delta_1,z_{n-1}^2-\Delta_{n-1}$ is of codimension $\ge1$.

By Lemmas~\ref{lem.slnd} and \ref{lem.sln}, $(\det-1,d)$ is a prime ideal of $K[\Mat_n]$. So we have an embedding $K[\Mat_n]/(\det-1,d)\to K[\mc V]$ which is the comorphism of a finite surjective morphism of varieties ${\mc V}\to V(\det-1,d)$, where $V(\det-1,d)$ is the closed subvariety of $\Mat_n$ that consists of the matrices of determinant 1 on which $d$ vanishes. This morphism maps the closed subvariety of $\mc V$ defined by the Jacobian ideal of $\det-1,d, f_1-s_1,\ldots,f_{n-1}-s_{n-1}, z_1^2-\Delta_1,\ldots,z_{n-1}^2-\Delta_{n-1}$ into the closed subvariety of $V(\det-1,d)$ defined by the ideal generated by the $2n$-th order minors of the Jacobian matrix of $(s_1,\ldots,s_n, d, \Delta_1,\ldots,\Delta_{n-1})$ with respect to the variables $\xi_{ij}$. This follows easily from the fact that $s_n=\det$ and that the $z_j$ and $u_j$ do not appear in the $s_i$ and $\Delta_i$. Since finite morphisms preserve dimension (see e.g. \cite{Eis} Corollary~9.3), it suffices to show that the latter variety is of codimension $\ge 1$ in $V(\det-1,d)$. Since $V(\det-1,d)$ is irreducible, this follows from Lemma~\ref{lem.jac}(ii).

\noindent 2.\ $\ov{B}_0(K)$ is an integral domain for any field $K$ of characteristic $p$.

We may assume that $K$ is algebraically closed. From the construction of the $f_i$ (see the proof of Lemma~\ref{lem.leadmonfi}) and the additivity of the $p$-th power map in characteristic $p$ it follows that $f_i=u_i^l$ mod $p$. So the kernel of the natural homomorphism from the polynomial algebra $K[(\xi_{ij})_{ij},u_1,\ldots,u_{n-1},z_1,\ldots,z_{n-1}]$ to $\ov{B}(K)$ is generated by the elements $\det-1,d,u_1^l-s_1,\ldots,u_{n-1}^l-s_{n-1}$ and the $\ov{A}(K)$-span of the monomials $\ov{u}_1^{k_1}\cdots \ov{u}_t^{k_t}$, $0\leq k_i<l$, is closed under multiplication for each $t\in\{0,\ldots,n-1\}$. We show by induction on $t$ that $\ov{B}_{0,t}(K):=\ov{A}(K)[\ov{u}_1,\ldots,\ov{u}_t]$ is an integral domain for $t=0,\ldots,n-1$. For $t=0$ this follows from Lemma~\ref{lem.slnd} and Proposition~\ref{prop.presentation(d)}(ii). Let $t\in\{1,\ldots,n-1\}$ and assume that it holds for $t-1$. Clearly $\ov{B}_{0,t}(K)=\ov{B}_{0,t-1}(K)[\ov{u}_t]\cong \ov{B}_{t-1}(K)[x]/(x^l-\ov{s}_t)$. So it suffices to prove that $x^l-\ov{s}_t$ is irreducible over the field of fractions of $\ov{B}_{0,t-1}(K)$. By the Vahlen-Capelli criterion or a more direct argument, it suffices to show that $\ov{s}_t$ is not a $p$-th power in the field of fractions of $\ov{B}_{0,t-1}(K)$. So assume that $\ov{s}_t=(v/w)^p$ for some $v,w\in\ov{B}_{0,t-1}(K)$ with $w\ne0$. Then we have $v^p=\ov{s}_tw^p=\ov{u}_t^lw^p$. So with $l'=l/p$, we have $(v-\ov{u}_t^{l'}w)^p=0$. But then $v-\ov{u}_t^{l'}w=0$ by Step 1. Now recall that $v$ and $w$ can be expressed uniquely as $\ov{A}(K)$-linear combinations of monomials in $\ov{u}_1,\ldots,\ov{u}_{t-1}$ with exponents $<l$. If such a monomial appears with a nonzero coefficient in $w$, then $\ov{u}_t^{l'}$ times this monomial appears with the same coefficient in the expression of $0=v-\ov{u}_t^{l'}w$ as an $\ov{A}(K)$-linear combination of restricted monomials in $\ov{u}_1,\ldots,\ov{u}_{n-1}$. Since this is impossible, we must have $w=0$. A contradiction.

\noindent 3.\ $\ov{B}_0({\mb C})$ is an integral domain.

This follows immediately from Step 2 and Lemma~\ref{lem.modp} applied to the $p$-adic valuation of $\mb Q$ and with $L=\mb C$.

\noindent 4.\ $\ov{B}_t({\mb C})$ is an integral domain for $t=0,\ldots,n-1$.

We prove this by induction on $t$. For $t=0$ it is the assertion of Step 3. Let $t\in\{1,\ldots,n-1\}$ and assume that it holds for $t-1$. Clearly $\ov{B}_t({\mb C})=\ov{B}_{t-1}({\mb C})[\ov{z}_t]\cong \ov{B}_{t-1}({\mb C})[x]/(x^2-\ov{\Delta}_t)$. So it suffices to prove that $x^2-\ov{\Delta}_t$ is irreducible over the field of fractions of $\ov{B}_{t-1}({\mb C})$. Assume that $x^2-\ov{\Delta}_t$ has a root in this field, i.e. that $\ov{\Delta}_t=(v/w)^2$ for some $v,w\in \ov{B}_{t-1}({\mb C})$ with $w\ne0$. By the same arguments as in the proof of Lemma~\ref{lem.algnum} we may assume that for some finite extension $F$ of $\mb Q$ there exist $v,w\in \ov{B}_{t-1}(F)$ with $w\ne0$ and $w^2\ov{\Delta}_t=v^2$. Let $\nu_2$ be an extension to $F$ of the $2$-adic valuation of $\mb Q$, let $S_2$ be the valuation ring of $\nu_2$, let $K$ be the residue class field and let $\delta\in S_2$ be a uniformiser for $\nu_2$. We may assume that the coefficients of $v$ and $w$ with respect to the ${\mb Z}$-basis of $\ov{B}_{t-1}$ mentioned earlier are in $S_2$. Assume that the coefficients of $w$ are all divisible by $\delta$ (in $S_2$). Then $w=0$ in $\ov{B}_{t-1}(K)$ and therefore $v^2=0$ in $\ov{B}_{t-1}(K)$. But by Step 1, $\ov{B}_{t-1}(K)$ is reduced, so $v=0$ in $\ov{B}_{t-1}(K)$ and all coefficients of $v$ are divisible by $\delta$. So, by cancelling a suitable power of $\delta$ in $w$ and $v$, we may assume that not all coefficients of $w$ are divisible by $\delta$. By passing to the residue class field $K$ we then obtain $v,w\in \ov{B}_{t-1}(K)$ with $w\ne0$ and $w^2\ov{\Delta}_t=v^2$. But then $(w\ov{z}_t-v)^2=0$ in $\ov{B}_t(K)$, since $\ov{z}_t^2=\ov{\Delta}_t$ and $K$ is of characteristic $2$. The reducedness of $\ov{B}_t(K)$ (Step 1) now gives $w\ov{z}_t-v=0$ in $\ov{B}_t(K)$. Now recall that $v$ and $w$ can be expressed uniquely as $\ov{A}(K)$-linear combinations of the monomials $\ov{u}_1^{k_1}\cdots \ov{u}_{n-1}^{k_{n-1}}\ov{z}_1^{m_1}\cdots \ov{z}_{t-1}^{m_{t-1}}$, $0\leq k_i<l$, $0\leq m_i<2$. We then obtain a contradiction in the same way as at the end of Step 2.

\noindent 5.\ $Z/(d)$ is an integral domain.

Since $\ov{Z}=\ov{B}({\mb C})[\ov{\Delta}^{\ -1}_1,\ldots,\ov{\Delta}^{\ -1}_{n-1}]$ and the $\ov{\Delta}_i$ are nonzero in $\ov{A}({\mb C})\cong {\mb C}[\SL_n]/(d')$ by Lemma~\ref{lem.slnd}, this follows from Step 4.
\end{proof}
\begin{remsnn}\ \\
1.\ Note that we didn't prove that $\ov{B}(K)$ is an integral domain for $K$ some algebraically closed  field of positive characteristic.\\
2.\ To attempt a proof for arbitrary odd $l>1$ I have tried the filtration with $\deg(\xi_{ij})=2l$, $\deg(z_i)=li$ and $\deg(u_i)=2i$. But the main problem with this filtration is that it does not simplify the relations $s_i=f_i(u_1,\ldots,u_{n-1})$ enough.
\end{remsnn}

 \end{document}